\documentclass[twoside,10pt]{article}
\usepackage{amsfonts}
\usepackage{amssymb}
\usepackage{amsmath}

\newtheorem{theorem}{Theorem}[section]

\newtheorem{lemma}[theorem]{Lemma}
\newtheorem{corollary}[theorem]{Corollary}
\newtheorem{definition}[theorem]{Definition}
\newtheorem{remark}[theorem]{Remark}
\newtheorem{proposition}[theorem]{Proposition}

\input{tcilatex}
\begin{document}

\title{Remarks on the $\Gamma $--regularization \\ of Non--convex and Non--semi--continuous \\  Functions on Topological Vector Spaces}
\author{J.-B. Bru and W. de Siqueira Pedra}
\date{\today }
\maketitle



\begin{abstract}
We show that the minimization problem of any non--convex and non--lower semi--continuous function on a compact convex subset
of a locally convex real topological vector space can be studied via
an associated convex and lower semi--continuous function $\Gamma \left( h\right) $. This observation uses the notion of $\Gamma $--regularization as a key ingredient.
As an application we obtain, on any locally convex real space, a generalization of the Lanford III--Robinson theorem which has
 only been proven for separable real Banach spaces. The latter is a characterization of subdifferentials of convex continuous functions. \\[0.3ex]
{\small \textit{Keywords:} variational problems, non--linear analysis, non--convexity, \\ $\Gamma $--regularization, Lanford III -- Robinson theorem.}\\[0.3ex]
\end{abstract}


\section{Introduction and Main Results}

\setcounter{equation}{0}%
Minimization problems $\inf \,h\left( K\right) $ on compact convex subsets $%
K $ of a locally convex real (topological vector) space\footnote{%
We assume throughout this paper that topological vector spaces are
Hausdorff spaces, i.e., points in those spaces define closed sets.}
 $\mathcal{X}$ are extensively studied for convex
and lower semi--continuous real--valued functions $h$. See, for instance,
\cite{Zeidler3}.

Such variational problems are, however, not systematically studied for \emph{%
non}--convex and \emph{non}--lower semi--continuous  real--valued functions $h$,
except for a few specific functions. See for instance \cite{Mueller}.
The aim of this paper is to show that
-- independently of convexity or lower semi--continuity of functions $h$
-- the minimization problem $\inf \,h\left( K\right) $ on compact convex
subsets $K$ of a locally convex real space $\mathcal{X}$ can be analyzed via
another minimization problem $\inf \,\Gamma \left( h\right) \left( K\right) $
associated with a convex and lower semi--continuous function $\Gamma
\left( h\right) $, for which various methods of analysis are available.

We are particularly interested in characterizing the following set of
generalized minimizers of any real--valued function $h$ on a compact convex set $K
$:

\begin{definition}[Set of generalized minimizers]
\label{gamm regularisation copy(4)}\mbox{ }\newline
Let $K$ be a (non--empty) compact convex subset of a locally convex real
space $\mathcal{X}$ and $h:K\rightarrow \left( -\infty ,\infty \right] $ be
any extended real--valued function. Then the set $\overline{\mathit{\Omega }\left(
h,K\right) }\subset K$ of generalized minimizers of $h$ is the closure of
the set
\begin{equation*}
\mathit{\Omega }\left( h,K\right) :=%
\Big\{%
x\in K:\exists \{x_{i}\}_{i\in I}\subset K\mathrm{\ \ with\ }%
x_{i}\rightarrow x\;\mathrm{and\;}\lim_{I}h(x_{i})=\inf \,h(K)%
\Big\}%
\end{equation*}%
of all limit points of approximating minimizers of $h$.
\end{definition}

\noindent Here, $\{x_{i}\}_{i\in I}\subset K$ is per definition a net of
\emph{approximating minimizers} when
\begin{equation*}
\underset{I}{\lim }\ h(x_{i})=\inf \,h(K).
\end{equation*}%
Note that, for any compact set $K$, $\mathit{\Omega }\left( h,K\right) $ is
non--empty because any net $\{x_{i}\}_{i\in I}\subset K$ converges along a
subnet.

In order to motivate the issue here, observe that $\inf \,h\left( K\right) $
can always be studied via a minimization problem associated with a (possibly
not convex, but) lower semi--continuous function $h_{0}$, known as the
\emph{lower semi--continuous hull} of $h$:

\begin{lemma}[Minimization of real--valued functions -- I]
\label{theorem trivial sympa 1 copy(2)}\mbox{ }\newline
Let $K$ be any (non--empty) compact, convex, and metrizable subset of a
locally convex real space $\mathcal{X}$ and $h:K\rightarrow \lbrack \mathrm{k%
},\infty ]$ be any extended real--valued
 function with $\mathrm{k}\in \mathbb{R}$. Then
there is a lower semi--continuous extended function
 $h_{0}:K\rightarrow \lbrack
\mathrm{k},\infty ]$ such that%
\begin{equation*}
\inf \,h\left( K\right) =\inf \,h_{0}\left( K\right) \quad \text{and}\quad
\mathit{\Omega }\left( h_{0},K\right) =\mathit{\Omega }\left( h,K\right) .
\end{equation*}
\end{lemma}

\noindent By lower semi--continuity, note that $\mathit{\Omega }\left(
h_{0},K\right) $ corresponds to the set of usual minimizers of $h_{0}$. Note
further that Lemma \ref{theorem trivial sympa 1 copy(2)} implies -- \ in the
case $K$ is metrizable -- that $\mathit{\Omega }\left( h,K\right) $ is
closed, again by lower semi--continuity of $h_{0}$. The proof of this lemma
is straightforward and is given in Section \ref{Section Proofs-I} for
completeness.

This result has two drawbacks: The compact convex set $K$ must be \emph{%
metrizable} in the elementary proof we give here and, more important, the
lower semi--continuous hull $h_{0}$ of $h$ is \emph{generally not convex}.
We give below a more elaborate result and show that both problems mentioned
above can be overcome by using the so--called $\Gamma $--regularization of
extended real--valued functions.
 The last is defined from the space $\mathrm{A}\left(
\mathcal{X}\right) $ of all affine continuous real--valued functions on a
locally convex real space\ $\mathcal{X}$ as follows (cf. \cite[Eq. (1.3) in
Chapter I]{Alfsen}):

\begin{definition}[$\Gamma $--regularization of real--valued functions]
\label{gamm regularisation}\mbox{ }\newline
For any extended real--valued function $h:K\rightarrow \lbrack \mathrm{k},\infty ]$ defined
on a (non--empty) compact convex subset $K\subset \mathcal{X}$, its $\Gamma $%
--regularization $\Gamma \left( h\right) $ on $K$ is the function defined
as the supremum over all affine and continuous minorants $m:\mathcal{X}%
\rightarrow \mathbb{R}$ of $h$, i.e., for all $x\in K$,
\begin{equation*}
\Gamma \left( h\right) \left( x\right) :=\sup \left\{ m(x):m\in \mathrm{A}%
\left( \mathcal{X}\right) \;\text{and }m|_{K}\leq h\right\} .
\end{equation*}
\end{definition}

\noindent Since the $\Gamma $--regularization $\Gamma \left( h\right) $ of a
extended
real--valued function $h$ is a supremum over continuous functions, $\Gamma
\left( h\right) $ is a \emph{convex} and \emph{lower semi--continuous}
function on $K$. For convenience, note that we identify extended
 real--valued functions $g $
 only defined on a\ convex compact subset $K\subset \mathcal{X}$ of the
locally convex real space\ $\mathcal{X}$ with its (trivial) extension
 $g_{\mathrm{ext}}$ to the whole space\ $\mathcal{X}$ defined by%
\begin{equation*}
g_{\mathrm{ext}}(x):=\left\{
\begin{array}{c}
g(x) \\
\infty%
\end{array}%
\begin{array}{l}
\mathrm{for}\ x\in K, \\
\mathrm{otherwise.}%
\end{array}%
\right.
\end{equation*}%
Clearly, with this prescription $g$ is lower semi--continuous (resp. convex) on $K$
iff $g$ is lower semi--continuous (resp. convex) on $\mathcal{X}$.

We prove in Section \ref{Section Proofs-II} the main result of this paper:

\begin{theorem}[Minimization of real--valued functions -- II]
\label{theorem trivial sympa 1}\mbox{ }\newline
Let $K$ be any (non--empty) compact convex subset of a locally convex real
space $\mathcal{X}$ and $h:K\rightarrow \lbrack \mathrm{k},\infty ]$ be any
extended real--valued function with $\mathrm{k}\in \mathbb{R}$. Then we have that:\newline
\emph{(i)}
\begin{equation*}
\inf \,h\left( K\right) =\inf \,\Gamma \left( h\right) \left( K\right) .
\end{equation*}%
\emph{(ii) }The set $\mathit{M}$ of minimizers of $\Gamma \left( h\right) $
over $K$ equals the closed convex hull of the set $\mathit{\Omega }\left(
h,K\right) $ of generalized minimizers of $h$ over $K$, i.e.,
\begin{equation*}
\mathit{M}=\overline{\mathrm{co}}\left( \mathit{\Omega }\left( h,K\right)
\right) .
\end{equation*}
\end{theorem}

This general fact related to the minimization of non--convex and non--lower
semi--continuous  real--valued
 functions on compact convex sets has not been
observed\footnote{Assertion (i) is, however, trivial.} before,
at least to our knowledge. Note that related results were obtained in
\cite{Benoist}\footnote{We thank the referee for pointing out this reference.}
for $\mathcal{X}=\mathbb{R}^n$. It
turns out to be extremely useful. It is, for instance, an essential argument
in the proof given in \cite{BruPedra2} of the validity of the so--called
Bogoliubov approximation on the level of states for a class of models for
fermions on the lattice. This problem, well--known in mathematical physics,
was first addressed by Ginibre \cite[p. 28]{Ginibre} in 1968 and is still
open for many physically important models.

Then, by using the theory of compact convex subsets of locally convex real
spaces $\mathcal{X}$ (see, e.g., \cite{Alfsen}), Theorem \ref{theorem
trivial sympa 1} yields a characterization of the set $\overline{\mathit{%
\Omega }\left( h,K\right) }$ of all generalized minimizers of $h$ over $K$.
Indeed, one important observation concerning locally convex real spaces $%
\mathcal{X}$ is that any compact convex subset $K\subset \mathcal{X}$ is the
closure of the convex hull of the (non--empty) set $\mathcal{E}(K)$ of its
extreme points, i.e., of the points which cannot be expressed as
(non--trivial) convex combinations of other elements in $K$. This is the
Krein--Milman theorem, see, e.g., \cite[Theorems 3.4 (b) and 3.23]{Rudin}.
In fact, among all subsets $Z\subset K$ generating $K$, $\mathcal{E}(K)$ is
-- in a sense -- the smallest one. This is the Milman theorem, see, e.g.,
\cite[Theorem 3.25]{Rudin}. It follows from Theorem \ref{theorem trivial
sympa 1} together with \cite[Theorems 3.4 (b), 3.23, 3.25]{Rudin} that
extreme points of the compact convex set
$\mathit{M}$ of minimizers of $\Gamma \left(h\right) $ over $K$
are generalized minimizers of $h$:

\begin{theorem}[Minimization of real-valued functions -- III]
\label{theorem trivial sympa 1 copy(1)}\mbox{ }\newline
Let $K$ be any (non--empty) compact convex subset of a locally convex real
space $\mathcal{X}$ and $h:K\rightarrow \lbrack \mathrm{k},\infty ]$ be any
extended real--valued
 function with $\mathrm{k}\in \mathbb{R}$. Then extreme points of the
compact convex set $\mathit{M}$ belong to the set of generalized minimizers
of $h$, i.e., $\mathcal{E}\left( \mathit{M}\right) \subseteq \overline{%
\mathit{\Omega }\left( h,K\right) }$.
\end{theorem}

\noindent This last result makes possible a \emph{full characterization} of
the closure of the set $\mathit{\Omega }\left( h,K\right) $ in the following
sense: Since $\mathit{M}$ is compact and convex, we can study the
minimization problem $\inf \,h\left( K_{\mathit{M}}\right) $ for any closed
(and hence compact) convex subset $K_{\mathit{M}}\subset \mathit{M}$.
Applying Theorem \ref{theorem trivial sympa 1} we get
\begin{equation}
\inf \,h\left( K_{\mathit{M}}\right) =\inf \,\Gamma \left( h|_{K_{\mathit{M}%
}}\right) \left( K_{\mathit{M}}\right) .  \label{full characterization1}
\end{equation}%
If
\begin{equation*}
\inf \,h\left( K_{\mathit{M}}\right) =\inf \,h\left( K\right)
\end{equation*}%
then, by Theorem \ref{theorem trivial sympa 1 copy(1)},
\begin{equation*}
\mathcal{E}\left( \mathit{M}_{K_{\mathit{M}}}\right) \subseteq \overline{%
\mathit{\Omega }\left( h|_{K_{\mathit{M}}},K_{\mathit{M}}\right) }\subseteq
\overline{\mathit{\Omega }\left( h,K\right) },
\end{equation*}%
where $\mathit{M}_{K_{\mathit{M}}}$ is the compact convex set of minimizers
of $\Gamma \left( h|_{K_{\mathit{M}}}\right) $ over $K_{\mathit{M}}\subset
\mathit{M}$. In general, $\mathcal{E}\left( \mathit{M}_{K_{\mathit{M}%
}}\right) \backslash \mathcal{E}\left( \mathit{M}\right) \neq \emptyset $
because $\mathit{M}_{K_{\mathit{M}}}$ is not necessarily a face of $\mathit{M%
}$. Thus we discover in this manner new points of $\overline{\mathit{\Omega }%
\left( h,K\right) }$ not contained in $\mathcal{E}\left( \mathit{M}\right) $%
. Choosing a sufficiently large family $\{K_{\mathit{M}}\}$ of closed convex
subsets of $\mathit{M}$ we can exhaust the set $\overline{\mathit{\Omega }%
\left( h,K\right) }$ through the union $\cup $ $\{\mathcal{E}\left( \mathit{M%
}_{K_{\mathit{M}}}\right) \}$. Note that this construction can be performed
in an inductive way: For each set $\mathit{M}_{K_{\mathit{M}}}$ of
minimizers consider further closed convex subsets $K_{\mathit{M}}^{\prime
}\subset $ $\mathit{M}_{K_{\mathit{M}}}$. The art consists in choosing the
family $\{K_{\mathit{M}}\}$ appropriately, i.e., it should be as small as
possible and the extreme points of $\mathit{M}_{K_{\mathit{M}}}$ should
possess some reasonable characterization. Of course, the latter   heavily
depends on the function $h$ and on particular properties of the compact
convex set $K$ (e.g., density of $\mathcal{E}(K)$, metrizability, etc.).

To close this section we recall that the $\Gamma $--regularization $\Gamma
\left( h\right) $ of a function $h$ on $K$ equals its twofold \emph{%
Legendre--Fenchel transform} -- also called the \emph{biconjugate }%
(function) of $h$. See, for instance,
\cite[Paragraph 51.3]{Zeidler3}. Indeed, $\Gamma \left( h\right) $ is the largest lower
semi--continuous and convex minorant of $h$ (cf. Corollary \ref{Biconjugate}%
). However, in contrast to the $\Gamma $--regularization the notion of
Legendre--Fenchel transform requires the use of dual pairs (cf. Definition %
\ref{dual pairs}). Since, for any locally convex real space $\mathcal{X}$
together with the space $\mathcal{X}^{\ast }$ of linear continuous
funtionals $\mathcal{X} \to \mathbb{R}$ (dual space)
equipped with the weak$^{\ast }$--topology,
$(\mathcal{X},\mathcal{X}^{\ast })$ is a dual pair, the
Legendre--Fenchel transform can be defined on any locally
convex real space $\mathcal{X}$ as follows:

\begin{definition}[The Legendre--Fenchel transform]
\label{Legendre--Fenchel transform}\mbox{ }\newline
Let $K$ be a (non--empty) compact convex subset of a locally convex real
space $\mathcal{X}$. For any extended real--valued
 function $h:K\rightarrow \left( -\infty,\infty \right] $, $h\not\equiv \infty$,
 its Legendre--Fenchel transform $h^{\ast }$ is the convex
weak$^{\ast }$--lower semi--continuous extented
function from $\mathcal{X}^{\ast }$
to $\left( -\infty ,\infty \right] $ defined, for any $x^{\ast }\in \mathcal{%
X}^{\ast }$, by
\begin{equation*}
h^{\ast }\left( x^{\ast }\right) :=\underset{x\in K}{\sup }\left\{ x^{\ast
}\left( x\right) -h\left( x\right) \right\} .
\end{equation*}
\end{definition}
See also \cite[Definition 51.1]{Zeidler3}.
Note that, together with its weak$^{\ast }$--topology, the dual
space $\mathcal{X}^{\ast }$ of any locally convex space $\mathcal{X}$ is
also a locally convex space, see \cite[Theorems 3.4 (b) and 3.10]{Rudin}.
Therefore, in case nothing is further specified, the space $\mathcal{X}%
^{\ast }$ is always equipped with its weak$^{\ast }$--topology.

The Legendre--Fenchel transform is strongly related to the notion of
Fenchel subdifferentials (see also \cite{Phelps-conv}):

\begin{definition}[Fenchel subdifferentials]
\label{tangent functional}\mbox{ }\newline
Let $h:\mathcal{X}\rightarrow \mathbb{(-\infty },\infty ]$ be any
extended real--valued
function on a real topological vector space $\mathcal{X}$. A
continuous linear
functional $\mathrm{d}h_{x}\in \mathcal{X}^{\ast }$ is said to be a
Fenchel subgradient (or tangent) of the function $h$ at $x\in \mathcal{X}$ iff,
for all $x^{\prime }\in \mathcal{X}$, $h(x+x^{\prime })\geq h(x)+\mathrm{d}%
h_{x}(x^{\prime })$. The set $\partial h(x)\subset \mathcal{X}^{\ast }$ of
Fenchel subgradients of $h$ at $x$ is called Fenchel
 subdifferential of $h$ at $x$.
\end{definition}

\noindent  Theorem \ref{theorem trivial sympa 1} establishes a
link between generalized minimizers and Fenchel subdifferentials:

\begin{theorem}[Subdifferentials of continuous convex functions -- I]
\label{theorem trivial sympa 3}Let $K$ be any (non--empty) compact convex
subset of a locally convex real space $\mathcal{X}$ and $h:K\rightarrow
\lbrack \mathrm{k},\infty ]$ be any extended real--valued
 function with $\mathrm{k}\in \mathbb{R}$. Then the Fenchel subdifferential
 $\partial h^{\ast }(x^{\ast })\subset \mathcal{X}$ of $h^{\ast }$
 at the point $x^{\ast }\in \mathcal{X}^{\ast }$
is the (non--empty) compact convex set
\begin{equation*}
\partial h^{\ast }(x^{\ast })=\overline{\mathrm{co}}\left( \mathit{\Omega }%
\left( h-x^{\ast },K\right) \right) .
\end{equation*}
\end{theorem}

\noindent This last result -- proven in Section \ref{Section Proofs-III} --
generalizes the Lanford III--Robinson theorem \cite[Theorem 1]{LanRob} which
has only been proven for separable real Banach spaces $\mathcal{X}$ and
continuous convex functions $h:\mathcal{X}\rightarrow \mathbb{R}$, cf.
Theorem \ref{Land.Rob}.

Indeed, for any extended real--valued function
 $h$ from a compact convex subset $K\subset \mathcal{X}$
 of a locally convex real space
$\mathcal{X}$ to $\left( -\infty ,\infty \right] $, let%
\begin{equation*}
\mathcal{Y}^{\ast }:=\left\{ x^{\ast }\in \mathcal{X}^{\ast }:h^{\ast }\text{
has a unique Fenchel subgradient }\mathrm{d}h_{x^{\ast }}^{\ast }\in \mathcal{X}%
\text{ at }x^{\ast }\right\} .
\end{equation*}%
For all $x^{\ast }\in \mathcal{X}^{\ast }$ and any open neighborhood $%
\mathcal{V}$ of\ $\{0\}\subset \mathcal{X}^{\ast }$, we also define the set%
\begin{equation}
\mathcal{T}_{x^{\ast },\mathcal{V}}:=\overline{\left\{ \mathrm{d}h_{y^{\ast
}}^{\ast }:y^{\ast }\in \mathcal{Y}^{\ast }\cap (x^{\ast }+\mathcal{V)}%
\right\} }^{\mathcal{X}}\subset \mathcal{X}  \label{set t voisinage}
\end{equation}%
and denote by $\mathcal{T}_{x^{\ast }}$ the intersection%
\begin{equation}
\mathcal{T}_{x^{\ast }}:=\bigcap\limits_{\mathcal{V}\ni 0\text{ open}}%
\mathcal{T}_{x^{\ast },\mathcal{V}}.  \label{set t voinagebis}
\end{equation}%
Here, $\overline{\;\cdot \;}^{\mathcal{X}}$ denotes the closure w.r.t. the
topology of $\mathcal{X}$.
Then we observe first that Theorem \ref{theorem trivial sympa 3} implies
that the set $\partial h^{\ast }(x^{\ast })\subset \mathcal{X}$ of Fenchel
subgradients of $h^{\ast }$ at the point $x^{\ast }\in \mathcal{X}^{\ast }$
is included in the closed convex hull of the set $\mathcal{T}_{x^{\ast }}$
provided $\mathcal{Y}^{\ast }$ is dense in $\mathcal{X}^{\ast }$ (cf.
Section \ref{Section Proofs-IV}):

\begin{corollary}[Subdifferentials of continuous convex functions -- II]
\label{corollary explosion lanford-robinson}Let $K$ be any (non--empty)
compact convex subset of a locally convex real space $\mathcal{X}$ and $%
h:K\rightarrow \lbrack \mathrm{k},\infty ]$ be any extended real--valued
 function with $\mathrm{k}\in \mathbb{R}$.
If $\mathcal{Y}^{\ast }$ is dense in $\mathcal{X}^{\ast }$
then, for any $x^{\ast }\in \mathcal{X}^{\ast }$,
\begin{equation*}
\partial h^{\ast }(x^{\ast })\subseteq \overline{\mathrm{co}}\left( \mathcal{T%
}_{x^{\ast }}\right) .
\end{equation*}
\end{corollary}

\noindent This last result applied on separable Banach spaces yields, in
turn, the following assertion (cf. Section \ref{Section Proofs-IV copy(1)}):

\begin{corollary}[The Lanford III--Robinson theorem]
\label{corollary explosion lanford-robinson copy(1)}\mbox{ }\newline
Let $\mathcal{X}$ be a separable Banach space and $h:\mathcal{X}\rightarrow
\mathbb{R}$ be any convex function which is globally Lipschitz continuous.
If the set%
\begin{equation*}
\mathcal{Y}:=\left\{ x\in \mathcal{X}:h\text{ has a unique Fenchel
subgradient }%
\mathrm{d}h_{x}\in \mathcal{X}^{\ast }\text{ at }x\right\}
\end{equation*}%
is dense in $\mathcal{X}$ then the Fenchel
subdifferential $\partial h(x)$ of $h$,
at any $x\in \mathcal{X}$, is the weak$^{\ast }$--closed convex hull of the
set $\mathcal{Z}_{x}$. Here, at fixed $x\in \mathcal{X}$, $\mathcal{Z}_{x}$
is the set of functionals $x^{\ast }\in \mathcal{X}^{\ast }$ such that there
is a net $\{x_{i}\}_{i\in I}$ in $\mathcal{Y}$ converging to $x$ with the
property that the unique Fenchel
subgradient $\mathrm{d}h_{x_{i}}\in \mathcal{X}^{\ast }$ of $h$ at
$x_{i}$ converges towards $x^{\ast }$ in the weak$^{\ast}$--topology.
\end{corollary}

\noindent Recall that the Mazur theorem shows that the set $\mathcal{Y}$ on
which a continuous convex function $h$ is G\^{a}teaux differentiable,
i.e., the set $\mathcal{Y}$ for which $h$ has exactly one Fenchel subgradient
\textrm{d}$h_{x}\in \mathcal{X}^{\ast }$ at any $x\in \mathcal{Y}$, is dense
in a separable Banach space $\mathcal{X}$, cf. Theorem \ref{Mazur} and
Remark \ref{Mazur remark}. Therefore, for globally Lipschitz continuous and
convex functions, the Lanford III--Robinson theorem \cite[Theorem 1]%
{LanRob} (cf. Theorem \ref{Land.Rob}) directly follows from Corollary \ref%
{corollary explosion lanford-robinson copy(1)}. Observe that, in which
concerns Fenchel
subdifferentials of convex continuous functions on Banach spaces,
the case of global Lipschitz continuous functions is already the most
general case: For any continuous convex function $h$ on a Banach space $%
\mathcal{X}$ and any $x\in \mathcal{X}$, there are $\varepsilon >0$ and a
globally Lipschitz continuous convex function $g$ such that $g\left(
y\right) =h\left( y\right) $ whenever $\left\Vert x-y\right\Vert
<\varepsilon $. In particular, $g$ and $h$ have the same Fenchel
subgradients at $x$.
Remark, indeed, that continuous convex functions $h$ on a Banach space
$\mathcal{X}$ are locally Lipschitz continuous and an example of such a
global Lipschitz continuous convex function is given by
\begin{equation*}
g\left( x\right) :=\inf \left\{ z\in \mathbb{R}:\left( x,z\right) \in \left[
\mathrm{epi}\left( h\right) +\mathcal{C}_{\alpha }\right] \right\} ,
\end{equation*}%
for sufficiently small $\alpha >0$. Here,
\begin{equation*}
\mathcal{C}_{\alpha }:=\left\{ \left( x,z\right) \in \mathcal{X}\times
\mathbb{R}:z\geq 0,\left\Vert x\right\Vert \leq \alpha z\right\}
\end{equation*}%
and $\mathrm{epi}\left( h\right) $ is the epigraph of $h$ defined by
\begin{equation*}
\mathrm{epi}\left( h\right) :=\left\{ \left( x,z\right) \in \mathcal{X}%
\times \mathbb{R}:z\geq f\left( x\right) \right\} .
\end{equation*}

The rest of the paper is structured as follows. Section \ref{Section Proofs}
gives the detailed proofs of Lemma \ref{theorem trivial sympa 1 copy(2)},
Theorems \ref{theorem trivial sympa 1}, \ref{theorem trivial sympa 3}, and
Corollaries \ref{corollary explosion lanford-robinson}--\ref{corollary
explosion lanford-robinson copy(1)}. Then, Section \ref{Concluding remarks}
discusses an additional observation which is relevant in the context of
minimization of non--convex or non--semi--continuous functions and which
does not seem to have been observed before. Indeed, Lemma \ref{Bauer maximum
principle bis} gives an extension of the Bauer maximum principle (Lemma \ref%
{Bauer maximum principle}). Finally, Section \ref{Section appendix} is a
concise appendix about dual pairs, barycenters
in relation with the $\Gamma $--regularization, the Mazur theorem, and the
Lanford III--Robinson theorem.

\section{Proofs\label{Section Proofs}}

This section gives the detailed proofs of Lemma \ref{theorem trivial sympa 1
copy(2)}, Theorems \ref{theorem trivial sympa 1}, \ref{theorem trivial sympa
3}, and Corollaries \ref{corollary explosion lanford-robinson}--\ref%
{corollary explosion lanford-robinson copy(1)}. Up to Corollary \ref%
{corollary explosion lanford-robinson copy(1)}, we will always assume that $%
K $ is a (non--empty) compact convex subset of a locally convex real space $%
\mathcal{X}$ and $h:K\rightarrow \lbrack \mathrm{k},\infty ]$ is any extended
real--valued
 function with $\mathrm{k}\in \mathbb{R}$. In Lemma \ref{theorem trivial
sympa 1 copy(2)} the metrizability of the topology on $K$ is also assumed.
In Corollary \ref{corollary explosion lanford-robinson copy(1)} $\mathcal{X}$
is a separable Banach space and $h:\mathcal{X}\rightarrow \mathbb{R}$ is any
globally Lipschitz continuous convex function.

\subsection{Proof of Lemma \protect\ref{theorem trivial sympa 1 copy(2)}
\label{Section Proofs-I}}

Because the subset $K\subset \mathcal{X}$ is metrizable and compact, it is
sequentially compact and we can restrict ourselves to sequences instead of
more general nets. Using any metric $d(x,y)$ on $K$ generating the topology
we define, at fixed $\delta >0$, the extended real--valued
function $h_{\delta }$ from $K$ to $[\mathrm{k},\infty ]$ by%
\begin{equation*}
h_{\delta }\left( x\right) :=\inf \,h(\mathcal{B}_{\delta }\left( x\right) )
\end{equation*}%
for any $x\in K$, where%
\begin{equation}
\mathcal{B}_{\delta }\left( x\right) :=\left\{ y\in K:\ d(x,y)<\delta
\right\}  \label{ball}
\end{equation}%
is the ball (in $K$) of radius $\delta >0$ centered at $x\in K$. The family $%
\{h_{\delta }\left( x\right) \}_{\delta >0}$ of extended real--valued
 functions is clearly
increasing as $\delta \searrow 0$ and is bounded from above by $h(x)$.
Therefore, for any $x\in K$, the limit of $h_{\delta }\left( x\right) \geq
\mathrm{k}$ as $\delta \searrow 0$ exists and defines an extended
 real--valued function
\begin{equation*}
x\mapsto h_{0}\left( x\right) :=\underset{\delta \searrow 0}{\lim }%
\,h_{\delta }\left( x\right)
\end{equation*}%
from $K$ to $[\mathrm{k},\infty ]$.

In fact, this construction is well--known and the function $h_{0}$ is
called the \emph{lower semi--continuous hull} of $h$ as it is a lower
semi--continuous extended real--valued function
 from $K$ to $[\mathrm{k},\infty ]$. Indeed, for
all $\delta >0$ and any sequence $\{x_{n}\}_{n=1}^{\infty }\subset K$
converging to $x\in K$, there is $N_{\delta }>0$ such that, for all $%
n>N_{\delta }$, $x_{n}\in \mathcal{B}_{\delta /2}\left( x\right) $ which
implies that $\mathcal{B}_{\delta /2}\left( x_{n}\right) \subset \mathcal{B}%
_{\delta }\left( x\right) $. In particular, $h_{\delta }\left( x\right) \leq
h_{\delta /2}\left( x_{n}\right) $ for all $\delta >0$ and $n>N_{\delta }$.
Since the family $\{h_{\delta }\left( x\right) \}_{\delta >0}$ defines an
increasing sequence as $\delta \searrow 0$, it follows that%
\begin{equation*}
h_{\delta }\left( x\right) \leq \ \liminf_{n\rightarrow \infty }h_{0}\left(
x_{n}\right)
\end{equation*}%
for any $\delta >0$ and $x\in K$. In the limit $\delta \searrow 0$ the
latter yields the lower semi--continuity of the extended real--valued
 function $h_{0}$ on $K $. Moreover,
\begin{equation}
h_{0}\left( x\right) \geq h_{\delta }\left( x\right) \geq \inf \,h(K)\geq
\mathrm{k}>-\infty  \label{definition de h lower3bis}
\end{equation}%
for any $x\in K$ and $\delta >0$.

We observe now that $h$ and $h_{0}$ have the same infimum on $K$:%
\begin{equation}
\inf h_{0}\left( K\right) =\inf h(x).  \label{definition de h lower3}
\end{equation}%
This can be seen by observing first that there is $y\in K$ such that
\begin{equation}
\inf h_{0}\left( K\right) =h_{0}\left( y\right)
\label{definition de h lower4}
\end{equation}%
because of the lower semi--continuity of $h_{0}$. Since $h_{\delta }\leq h$
on $K$ for any $\delta >0$, we have $h_{0}\leq h\ $on $K$, which combined
with (\ref{definition de h lower3bis}) and (\ref{definition de h lower4})
yields Equality (\ref{definition de h lower3}).

Additionally, for all $\delta >0$ and any minimizer $y\in K$ of $h_{0}$ over
$K$, there is a sequence $\{x_{\delta ,n}\}_{n=1}^{\infty }\subset \mathcal{B%
}_{\delta }\left( y\right) $ of approximating minimizers of $h$ over $%
\mathcal{B}_{\delta }\left( y\right) $, that is,
\begin{equation*}
h_{\delta }\left( y\right) :=\inf \,h(\mathcal{B}_{\delta }\left( y\right) )=%
\underset{n\rightarrow \infty }{\lim }h(x_{\delta ,n})\leq h(y).
\end{equation*}%
We can assume without loss of generality that
\begin{equation*}
d(x_{\delta ,n},y)\leq \delta \mathrm{\quad and\quad }|h(x_{\delta
,n})-h_{\delta }\left( y\right) |\leq 2^{-n}
\end{equation*}%
for all $n\in \mathbb{N}$ and all $\delta >0$. Note that $h_{\delta }\left(
y\right) \rightarrow $ $h_{0}\left( y\right) $ as $\delta \searrow 0$. Thus,
by taking any function $p(\delta )\in \mathbb{N}$ satisfying $p(\delta
)>\delta ^{-1}$ we obtain that $x_{\delta ,p(\delta )}$ converges to $y\in K$
as $\delta \searrow 0$ with the property that $h(x_{\delta ,p(\delta )})$
converges to $h_{0}\left( y\right) $. Using Equalities (\ref{definition de h
lower3}) and (\ref{definition de h lower4}) we obtain that all minimizers of
(\ref{definition de h lower4}) are generalized minimizers of $h$, i.e.,
\begin{equation*}
\mathit{\Omega }\left( h_{0},K\right) \subseteq \mathit{\Omega }\left(
h,K\right) .
\end{equation*}%
The converse inclusion%
\begin{equation*}
\mathit{\Omega }\left( h,K\right) \subseteq \mathit{\Omega }\left(
h_{0},K\right)
\end{equation*}%
is straightforward because one has the inequality $h_{0}\leq h$ on $K$ as
well as Equality (\ref{definition de h lower3}).

\subsection{Proof of Theorem \protect\ref{theorem trivial sympa 1}\label%
{Section Proofs-II}}

The assertion (i) of Theorem \ref{theorem trivial sympa 1} is a standard
result. Indeed, by Definition \ref{gamm regularisation}, $\Gamma \left(
h\right) \leq h$ on $K$ and thus%
\begin{equation*}
\inf \,\Gamma \left( h\right) \left( K\right) \leq \inf \,h\left( K\right) .
\end{equation*}%
The converse inequality is derived by restricting the supremum in Definition %
\ref{gamm regularisation} to constant maps $m$ from $\mathcal{X}$ to $%
\mathbb{R}$ with $\mathrm{k}\leq m\leq h$ on $K$.

Observe that the
variational problem
$\inf \,\Gamma \left( h\right) (K)$ has minimizers and
the set $\mathit{M}=\mathit{\Omega }\left( \Gamma \left( h\right) ,K\right) $
of all minimizers of $\Gamma \left( h\right) $ is convex and
compact.
For any $y\in \mathit{\Omega }\left( h,K\right) $, there is a net $\left\{
x_{i}\right\} _{i\in I}\subset K$ of approximating minimizers of $h$ on $K$
converging to $y$. In particular, since the function $\Gamma \left(
h\right) $ is lower semi--continuous and $\Gamma \left( h\right) \leq h$ on $%
K$, we have that
\begin{equation*}
\Gamma \left( h\right) (y)\leq \underset{I}{\liminf }\,\Gamma \left(
h\right) (x_{i})\leq \underset{I}{\lim }\,h(x_{i})=\inf \,h(K)=\inf \,\Gamma
\left( h\right) (K),
\end{equation*}%
i.e., $y\in \mathit{M}$. Since $\mathit{M}$ is convex and compact, we obtain
that
\begin{equation}
\mathit{M}\supset \overline{\mathrm{co}}\left( \mathit{\Omega }\left(
h,K\right) \right) .  \label{inclusion1}
\end{equation}%
So, we prove now the converse inclusion. We can assume without loss of
generality that $\overline{\mathrm{co}}\left( \mathit{\Omega }\left(
h,K\right) \right) \neq K$ since otherwise there is nothing to prove. We
show\ next that, for any $x\in K\backslash \overline{\mathrm{co}}\left(
\mathit{\Omega }\left( h,K\right) \right) $, we have $x\notin $ $\mathit{M}$%
.

As $\overline{\mathrm{co}}\left( \mathit{\Omega }\left( h,K\right) \right) $
is a closed set of a locally convex real space $\mathcal{X}$, for any $x\in
K\backslash \overline{\mathrm{co}}\left( \mathit{\Omega }\left( h,K\right)
\right)$, there is an open and convex neighborhood $\mathcal{V}_{x}\subset
$ $\mathcal{X}$ of $\{0\}\subset \mathcal{X}$ which is symmetric, i.e., $%
\mathcal{V}_{x}=-\mathcal{V}_{x}$, and which satisfies
\begin{equation*}
\mathcal{G}_{x}\cap \left[ \{x\}+\mathcal{V}_{x}\right] =\emptyset
\end{equation*}%
with%
\begin{equation*}
\mathcal{G}_{x}:=K\cap \left[ \overline{\mathrm{co}}\left( \mathit{\Omega }%
\left( h,K\right) \right) +\mathcal{V}_{x}\right] .
\end{equation*}%
This follows from \cite[Theorem 1.10]{Rudin} together with the fact that
each neighborhood of $\{0\}\subset \mathcal{X}$ contains some open and
convex neighborhood of $\{0\}\subset \mathcal{X}$ because $\mathcal{X}$ is
locally convex. Observe also that any one--point set $\{x\}\subset $ $%
\mathcal{X}$ is trivially compact.

For any neighborhood $\mathcal{V}_{x}$ of $\{0\}\subset \mathcal{X}$ in a
locally convex real space, there is another convex, symmetric, and open
neighborhood $\mathcal{V}_{x}^{\prime }$ of $\{0\}\subset \mathcal{X}$ such
that $[\mathcal{V}_{x}^{\prime }+\mathcal{V}_{x}^{\prime }]\subset \mathcal{V%
}_{x}$, see proof of \cite[Theorem 1.10]{Rudin}. Let%
\begin{equation*}
\mathcal{G}_{x}^{\prime }:=K\cap \left[ \overline{\mathrm{co}}\left( \mathit{%
\Omega }\left( h,K\right) \right) +\mathcal{V}_{x}^{\prime }\right] .
\end{equation*}%
Then the following inclusions hold:%
\begin{equation}
\overline{\mathrm{co}}\left( \mathit{\Omega }\left( h,K\right) \right) %
\subset \mathcal{G}_{x}^{\prime }\subset \overline{\mathcal{G}_{x}^{\prime }}%
\subset \mathcal{G}_{x}\subset \overline{\mathcal{G}_{x}}\subset K\backslash
\{x\}.  \label{eq sup}
\end{equation}%
Since $K$, $\mathcal{V}_{x}$, $\mathcal{V}_{x}^{\prime }$, and $\overline{%
\mathrm{co}} \left( \mathit{\Omega }\left( h,K\right) \right) $ are all
convex sets, $\mathcal{G}_{x}$ and $\mathcal{G}_{x}^{\prime }$ are also
convex. Seen as subsets of $K$ they are open neighborhoods of $\overline{%
\mathrm{co}}\left( \mathit{\Omega }\left( h,K\right) \right) $.

The set $\mathcal{X}$ is a
Hausdorff space and thus any compact subset $K$ of $\mathcal{X}$ is
a normal space. By Urysohn lemma, there is a continuous function%
\begin{equation*}
f_{x}:K\rightarrow \lbrack \inf h(K),\inf h(K\backslash \mathcal{G}%
_{x}^{\prime })]
\end{equation*}%
satisfying $f_{x}\leq h$ and
\begin{equation*}
f_{x}\left( y\right) =\left\{
\begin{array}{ll}
\inf h(K) & \mathrm{for\ }y\in \overline{\mathcal{G}_{x}^{\prime }}. \\
\inf h(K\backslash \mathcal{G}_{x}^{\prime }) & \mathrm{for\ }y\in
K\backslash \mathcal{G}_{x}.%
\end{array}%
\right.
\end{equation*}%
By compactness of $K\backslash \mathcal{G}_{x}^{\prime }$ and the inclusion
$\mathit{\Omega }\left( h,K\right) \subset \mathcal{G}_{x}^{\prime }$,
observe that%
\begin{equation*}
\inf h(K\backslash \mathcal{G}_{x}^{\prime })>\inf h(K).
\end{equation*}%
Then we have per construction that
\begin{equation}
f_{x}(\overline{\mathrm{co}}\left( \mathit{\Omega }\left( h,K\right) \right) %
)=\{\inf h(K)\}  \label{Omega f sympa}
\end{equation}%
and%
\begin{equation}
f_{x}^{-1}(\inf h(K))=\mathit{\Omega }\left( f_{x},K\right) \subset \mathcal{%
G}_{x}  \label{Omega f sympabis}
\end{equation}%
for any $x\in K\backslash \overline{\mathrm{co}}\left( \mathit{\Omega }\left(
h,K\right) \right) $.

We use now the $\Gamma $--regularization $\Gamma \left( f_{x}\right) $ of $%
f_{x}$ on the set $K$ and denote by $\mathit{M}_{x}=\mathit{\Omega }\left(
\Gamma (f_{x}),K\right) $ its non--empty set of minimizers over $K$.
Applying Theorem \ref{Thm - Corollary I.3.6}, for any $y\in \mathit{M}_{x}$,
we have a probability measure $\mu _{y}\in M_{1}^{+}(K)$\ on $K$ with
barycenter $y$ such that
\begin{equation}
\Gamma \left( f_{x}\right) \left( y\right) =\int_{K}\mathrm{d}\mu
_{y}(z)\;f_{x}\left( z\right) .  \label{herve bis}
\end{equation}%
As $y\in \mathit{M}_{x}$, i.e.,
\begin{equation}
\Gamma \left( f_{x}\right) \left( y\right) =\inf \,\Gamma \left(
f_{x}\right) (K)=\inf f_{x}(K),  \label{herve 2}
\end{equation}%
we deduce from (\ref{herve bis}) that
\begin{equation*}
\mu _{y}(\mathit{\Omega }\left( f_{x},K\right) )=1
\end{equation*}%
and it follows that $y\in $ $\overline{\mathrm{co}}\left( \mathit{\Omega }%
\left( f_{x},K\right) \right) $, by Theorem \ref{thm barycenter}. Using (%
\ref{Omega f sympabis}) together with the convexity of the open neighborhood
$\mathcal{G}_{x}$ of $\overline{\mathrm{co}}\left( \mathit{\Omega }\left(
h,K\right) \right) $ we thus obtain%
\begin{equation}
\mathit{M}_{x}\subset \overline{\mathrm{co}}\left( \mathit{\Omega }\left(
f_{x},K\right) \right) \subset \overline{\mathcal{G}_{x}}  \label{herve 3}
\end{equation}%
for any $x\in K\backslash \overline{\mathrm{co}}\left( \mathit{\Omega }\left(
h,K\right) \right) $.

We remark now that the inequality $f_{x}\leq h$ on $K$ yields $\Gamma \left(
f_{x}\right) \leq \Gamma \left( h\right) $ on $K$ because of Corollary \ref%
{Biconjugate}. As a consequence, it results from (i) and (\ref{Omega f sympa}%
) that the set $\mathit{M}$ of minimizers of $\Gamma \left( h\right) $ over $%
K$ is included in $\mathit{M}_{x}$, i.e., $\mathit{M}\subset \mathit{M}_{x}$%
. Hence, by (\ref{eq sup}) and (\ref{herve 3}), we have the inclusions
\begin{equation}
\mathit{M}\subset \overline{\mathcal{G}_{x}}\subset K\backslash \{x\}.
\label{inclusion2bis}
\end{equation}%
Therefore, we combine (\ref{inclusion1}) with (\ref{inclusion2bis}) for all $%
x\in K\backslash \overline{\mathrm{co}}\left( \mathit{\Omega }\left(
h,K\right) \right) $ to obtain the desired equality in the assertion (ii)
of Theorem \ref{theorem trivial sympa 1}.

\subsection{Proof of Theorem \protect\ref{theorem trivial sympa 3}\label%
{Section Proofs-III}}

The proof of Theorem \ref{theorem trivial sympa 3} is a simple consequence
of Theorem \ref{theorem trivial sympa 1} together with the following
well--known result:

\begin{lemma}[Fenchel subgradients as minimizers]
\label{theorem trivial sympa 2}\mbox{ }\newline
Let $(\mathcal{X},\mathcal{X}^{\ast })$ be a dual pair and $h\not\equiv \infty$
 be any extended real--valued
function from a (non--empty) convex subset $K\subseteq \mathcal{X}$ to $%
(-\infty ,\infty ]$. Then the Fenchel
subdifferential $\partial h^{\ast }(x^{\ast
})\subset \mathcal{X}$ of $h^{\ast }$ at the point $x^{\ast }\in \mathcal{X}%
^{\ast }$ is the (non--empty) set $\mathit{M}_{x^{\ast }}$ of minimizers
over $K$ of the map
\begin{equation*}
y\mapsto \Gamma \left( h\right) \left( y\right) -x^{\ast }\left( y\right)
\end{equation*}%
from $K\subseteq \mathcal{X}$ to $(-\infty ,\infty ]$.
\end{lemma}

\textit{Proof. }The proof is standard and simple, see, e.g., \cite[Theorem
I.6.6]{Simon}. Indeed, any Fenchel subgradient $x\in \mathcal{X}$ of the
Legendre--Fenchel transform $h^{\ast }$ at the point $x^{\ast }\in \mathcal{X%
}$ satisfies the inequality:
\begin{equation}
x^{\ast }\left( x\right) +h^{\ast }\left( y^{\ast }\right) -y^{\ast }\left(
x\right) \geq h^{\ast }\left( x^{\ast }\right)  \label{landford1landford1}
\end{equation}%
for any $y^{\ast }\in \mathcal{X}^{\ast }$, see Definition \ref{tangent
functional}. Since $h^{\ast }=h^{\ast \ast \ast }$ and $\Gamma \left(
h\right) =h^{\ast \ast }$ (cf. Corollary \ref{Biconjugate} and \cite[%
Proposition 51.6]{Zeidler3}), we have (\ref{landford1landford1}) iff
\begin{equation*}
x^{\ast }\left( x\right) +\underset{y^{\ast }\in \mathcal{X}^{\ast }}{\inf }%
\left\{ h^{\ast }\left( y^{\ast }\right) -y^{\ast }\left( x\right) \right\}
=x^{\ast }\left( x\right) -\Gamma \left( h\right) \left( x\right) \geq
\underset{y\in K}{\sup }\left\{ x^{\ast }\left( y\right) -\Gamma \left(
h\right) \left( y\right) \right\} ,
\end{equation*}%
see Definition \ref{Legendre--Fenchel transform}.
\hfill\qed\medbreak%

We combine now Theorem \ref{theorem trivial sympa 1} with Lemma \ref{theorem
trivial sympa 2} to characterize the Fenchel subdifferential $\partial h^{\ast
}(x^{\ast })\subset \mathcal{X}$ of $h^{\ast }$ at the point $x^{\ast }\in
\mathcal{X}^{\ast }$ as the closed convex hull of the set $\mathit{\Omega }%
\left( h-x^{\ast },K\right) $ of generalized minimizers of $h$ over a
compact convex subset $K$, see Definition \ref{gamm regularisation copy(4)}.
Indeed, for any $x^{\ast }\in \mathcal{X}^{\ast }$,
\begin{equation*}
\Gamma \left( h-x^{\ast }\right) =\Gamma \left( h\right) -x^{\ast },
\end{equation*}%
see Definition \ref{gamm regularisation}.

\subsection{Proof of Corollary \protect\ref{corollary explosion
lanford-robinson}\label{Section Proofs-IV}}

For $x^{\ast }\in \mathcal{X}^{\ast }$ and any open neighborhood $\mathcal{V}
$ of $\{0\}\subset \mathcal{X}^{\ast }$, we define the map $g_{\mathcal{V}%
,x^{\ast }}$ from $\mathcal{X}$ to $[\mathrm{k},\infty ]$ with $\mathrm{k}%
\in \mathbb{R}$ by%
\begin{equation*}
g_{\mathcal{V},x^{\ast }}\left( x\right) :=\left\{
\begin{array}{c}
\Gamma \left( h\right) \left( x\right) \\
\infty%
\end{array}%
\begin{array}{l}
\mathrm{for}\ x=\mathrm{d}h_{y^{\ast }}^{\ast }\mathrm{\ with}\ y^{\ast }\in
\mathcal{Y}^{\ast }\cap (x^{\ast }+\mathcal{V)}, \\
\mathrm{otherwise.}%
\end{array}%
\right.
\end{equation*}%
For any $y^{\ast }\in \mathcal{Y}^{\ast }\cap \left( x^{\ast }+\mathcal{V}%
\right) $, one has the equality $g_{\mathcal{V},x^{\ast }}^{\ast }\left(
y^{\ast }\right) =h^{\ast }\left( y^{\ast }\right) $. This easily follows
from the fact that
\begin{align*}
h^{\ast }\left( y^{\ast }\right) & =\underset{z\in K}{\sup }\left\{ y^{\ast
}\left( z\right) -\Gamma \left( h\right) \left( z\right) \right\} =y^{\ast
}\left( x\right) -\Gamma \left( h\right) \left( x\right) \\
& =\underset{z\in K}{\sup }\left\{ y^{\ast }\left( z\right) -g_{\mathcal{V}%
,x^{\ast }}\left( z\right) \right\} =g_{\mathcal{V},x^{\ast }}^{\ast }\left(
y^{\ast }\right)
\end{align*}%
with $x:=\mathrm{d}h_{y^{\ast }}^{\ast }$, see proof of Lemma \ref{theorem
trivial sympa 2}. Let $\mathcal{W}$ be any open neighborhood of $%
\{0\}\subset \mathcal{X}^{\ast }$. Then, for any $z\in K$, the set
\begin{equation*}
\{\delta ^{\ast }(z):\delta ^{\ast }\in \mathcal{W}\}\subset \mathbb{R}
\end{equation*}%
is bounded, by continuity of the linear map $\delta ^{\ast }\mapsto \delta
^{\ast }(z)$. From the the principle of uniform boundedness for compact
convex sets, i.e., the version of the Banach--Steinhaus theorem stated, for
instance, in \cite[Theorem 2.9]{Rudin}, the set%
\begin{equation*}
\{\delta ^{\ast }(z):\delta ^{\ast }\in \mathcal{W},\,z\in K\}\subset
\mathbb{R}
\end{equation*}%
is also bounded. Thus, for any $z^{\ast }\in \mathcal{X}^{\ast }$,%
\begin{eqnarray*}
\lim_{s\searrow 0}\sup \left\{ |h^{\ast }\left( z^{\ast }\right) -h^{\ast
}\left( z^{\ast }+\delta ^{\ast }\right) |\,:\,\delta ^{\ast }\in s\mathcal{W%
}\right\} &=&0, \\
\lim_{s\searrow 0}\sup \left\{ |g_{\mathcal{V},x^{\ast }}^{\ast }\left(
z^{\ast }\right) -g_{\mathcal{V},x^{\ast }}^{\ast }\left( z^{\ast }+\delta
^{\ast }\right) |\,:\,\delta ^{\ast }\in s\mathcal{W}\right\} &=&0.
\end{eqnarray*}%
This implies the continuity of $h^{\ast }$ and $g_{\mathcal{V},x^{\ast
}}^{\ast }$. Hence, from the density of $\mathcal{Y}^{\ast }$, $h^{\ast }=g_{%
\mathcal{V},x^{\ast }}^{\ast }$ on the open neighborhood $\left( x^{\ast }+%
\mathcal{V}\right) $ of $\{x^{\ast }\}\subset \mathcal{X}^{\ast }$. In
particular, $h^{\ast }$ and $g_{\mathcal{V},x^{\ast }}^{\ast }$ have the
same Fenchel subgradients at the point $x^{\ast }$. From Theorems \ref{theorem
trivial sympa 1 copy(1)} and \ref{theorem trivial sympa 3}, for each open
neighborhood $\mathcal{V}$ of $\{0\}\subset \mathcal{X}^{\ast }$, the
extreme Fenchel
subgradients of $h^{\ast }$ at $x^{\ast }$ are all contained in the
set $\mathcal{T}_{x^{\ast },\mathcal{V}}$ defined by (\ref{set t voisinage}%
). Corollary \ref{corollary explosion lanford-robinson} thus follows.

\subsection{Proof of Corollary \protect\ref{corollary explosion
lanford-robinson copy(1)}\label{Section Proofs-IV copy(1)}}

Note that $h^{\ast \ast }=h$ because the function $h$ is continuous and
convex. By the global Lipschitz continuity of $h$,
\begin{equation*}
h\left( x\right) =\underset{x^{\ast }\in \mathcal{X}^{\ast }}{\sup }\left\{
x^{\ast }\left( x\right) -h^{\ast }\left( x^{\ast }\right) \right\} =%
\underset{x^{\ast }\in K}{\sup }\left\{ x^{\ast }\left( x\right) -h^{\ast
}\left( x^{\ast }\right) \right\}
\end{equation*}%
with $K:=\mathcal{B}_{R}\left( 0\right) \subset \mathcal{X}^{\ast }$ being
some ball of sufficiently large radius $R>0$ centered at $0$. The set $K$ is
weak$^{\ast }$--compact, by the Banach--Alaoglu theorem.

Now, for any fixed $x\in \mathcal{X}$ and all $x^{\ast }\in \mathcal{Z}%
_{x}\subset \mathcal{X}^{\ast }$, by definition of the set $\mathcal{Z}_{x}$%
, there is a net $\{x_{i}\}_{i\in I}$ in $\mathcal{Y}$ converging to $x$
with the property that the unique Fenchel
subgradient $x_{i}^{\ast }:=\mathrm{d}h_{x_{i}}\in \mathcal{X}^{\ast }$
of $h$ at $x_{i}$ converges towards $x^{\ast }$
in the weak$^{\ast }$--topology. Therefore, by continuity of $h$,
for any fixed $x\in \mathcal{X}$ and all $x^{\ast }\in \mathcal{Z}_{x}$,%
\begin{equation*}
h\left( x\right) =\underset{y^{\ast }\in \mathcal{X}^{\ast }}{\sup }\left\{
y^{\ast }\left( x\right) -h^{\ast }\left( y^{\ast }\right) \right\} =%
\underset{I}{\lim }\ h\left( x_{i}\right) =\underset{I}{\lim }\left\{
x_{i}^{\ast }\left( x\right) -h^{\ast }\left( x_{i}^{\ast }\right) \right\} ,
\end{equation*}%
with $\{x_{i}^{\ast }\}_{i\in I}$ converging to $x^{\ast }$. In other words,%
\begin{equation*}
\mathcal{Z}_{x}\subset \mathit{\Omega }\left( h^{\ast }-x,K\right) ,
\end{equation*}%
see Definition \ref{gamm regularisation copy(4)}. Thus, by Theorem \ref%
{theorem trivial sympa 3} and Corollary \ref{corollary explosion
lanford-robinson}, it suffices to prove that $\mathcal{T}_{x}\subset
\mathcal{Z}_{x}$.

By density of $\mathcal{Y}$ in $\mathcal{X}$, observe that the set
\begin{equation*}
\mathcal{T}_{x,\mathcal{V}}:=\overline{\left\{ \mathrm{d}h_{y}:y\in \mathcal{%
Y}\cap (x+\mathcal{V)}\right\} }^{\mathcal{X}^{\ast }}\subset \mathcal{X}^{\ast }
\end{equation*}%
is non--empty for any open neighborhood $\mathcal{V}$ of $\{0\}\subset
\mathcal{X}$. Meanwhile, the weak$^{\ast }$--compact set $K$ is metrizable
with respect to (w.r.t.) the weak$^{\ast }$--topology, by separability of $%
\mathcal{X}$, see \cite[Theorem 3.16]{Rudin}. In particular, $K$ is
sequentially compact and we can restrict ourselves to sequences instead of
more general nets. In particular, by (\ref{set t voisinage})--(\ref{set t
voinagebis}), one has
\begin{equation}
\mathcal{T}_{x}=\bigcap\limits_{n\in \mathbb{N}}\mathcal{T}_{x,\mathcal{B}%
_{1/n}\left( 0\right) }  \label{inclusion sup}
\end{equation}%
with $\mathcal{B}_{\delta }\left( x\right) $ being the ball (in $K$) of
radius $\delta >0$ centered at $x\in K$. Here, $\mathcal{B}_{\delta }\left(
x\right) $ is defined by (\ref{ball}) for any metric $d$ on $K$ generating
its weak$^{\ast }$--topology. For any $x^{\ast }\in \mathcal{T}_{x}\subset K$
and any $n\in \mathbb{N}$, there are per definition a sequence $%
\{x_{n,m}^{\ast }\}_{m=1}^{\infty }$ converging to $x^{\ast }$ in $K$ as $%
m\rightarrow \infty $ and an integer $N_{n}>0$ such that, for all $m\geq
N_{n}$, $d(x^{\ast },x_{n,m}^{\ast })\leq 2^{-n}$ and $x_{n,m}^{\ast }=%
\mathrm{d}h_{x_{n,m}}$ for some $x_{n,m}\in \mathcal{Y}\cap \lbrack x+%
\mathcal{B}_{1/n}\left( 0\right) ]$. Taking any function $p(n)\in \mathbb{N}$
satisfying $p(n)>N_{n}$ and converging to $\infty $ as $n\rightarrow \infty $
we obtain a sequence $\{x_{n,p(n)}^{\ast }\}_{n=1}^{\infty }$ converging to $%
x^{\ast }\in \mathcal{Z}_{x}$ as $n\rightarrow \infty $. This yields the
inclusion $\mathcal{T}_{x}\subset \mathcal{Z}_{x}$.

\section{Further Remarks\label{Concluding remarks}}

We give here an additional observation which is\ not necessarily directly
related to the main results of the paper. It concerns an extension of the
Bauer maximum principle \cite[Theorem I.5.3.]{Alfsen}. See \cite{BruPedra2}
for an application to statistical mechanics.

First, recall that the $\Gamma $--regularization $\Gamma \left( h\right) $
of an extended real--valued function $h$ is a convex and lower semi--continuous function
on a compact convex subset $K$. Moreover, every convex and lower
semi--continuous function on $K$ equals its own $\Gamma $--regularization
on $K$ (see, e.g., \cite[Proposition I.1.2.]{Alfsen}):

\begin{proposition}[$\Gamma $--regularization of lower semi--cont. conv. maps%
]
\label{lemma gamma regularisation}Let $h$ be any extended function from a
(non--empty) compact convex subset $K\subset \mathcal{X}$ of a locally
convex real space $\mathcal{X}$ to $\left( -\infty ,\infty \right] $. Then
the following statements are equivalent:\newline
\emph{(i)} $\Gamma \left( h\right) =h$ on $K$.\newline
\emph{(ii)} $h$ is a lower semi--continuous convex function on $K$.
\end{proposition}

\noindent This proposition is a standard result. The compactness of $K$ is
in fact not necessary but $K$ should be a closed convex set. This result can
directly be proven without using the fact that the $\Gamma $--regularization
$\Gamma \left( h\right) $ of a function $h$ on $K$ equals its twofold
\emph{Legendre--Fenchel transform} -- also called the \emph{biconjugate }%
(function) of $h$. Indeed, $\Gamma \left( h\right) $ is the largest lower
semi--continuous and convex minorant of $h$:

\begin{corollary}[Largest lower semi--cont. convex minorant of $h$]
\label{Biconjugate}\mbox{ }\newline
Let $h$ be any extended function from a (non--empty) compact convex subset $%
K\subset \mathcal{X}$ of a locally convex real space $\mathcal{X}$ to $%
\left( -\infty ,\infty \right] $. Then its $\Gamma $--regularization $\Gamma
\left( h\right) $ is its largest lower semi--continuous convex minorant on $%
K $.
\end{corollary}

\textit{Proof. }For any lower semi--continuous convex extended real--valued
 function $f$
defined on $K$ satisfying $f\leq h$, we have, by Proposition \ref{lemma
gamma regularisation}, that
\begin{equation*}
f\left( x\right) =\sup \left\{ m(x):m\in \mathrm{A}\left( \mathcal{X}\right)
\;\text{and }m|_{K}\leq f\leq h\right\} \leq \Gamma \left( h\right) \left(
x\right)
\end{equation*}%
for any $x\in K$.
\hfill\qed\medbreak%

\noindent In particular, if $(\mathcal{X},\mathcal{X}^{\ast })$ is a dual
pair and $h \not\equiv \infty$ is any
extended function from $K$ to $(-\infty ,\infty ]$ then $\Gamma
\left( h\right) =h^{\ast \ast }$, see \cite[Proposition 51.6]{Zeidler3}.

Proposition \ref{lemma gamma regularisation} has another interesting
consequence: An extension of the Bauer maximum principle \cite[Theorem I.5.3.%
]{Alfsen} which, in the case of convex functions, is:

\begin{lemma}[Bauer maximum principle]
\label{Bauer maximum principle}\mbox{ }\newline
Let $\mathcal{X}$ be a locally convex real space. An upper semi--continuous
convex  real--valued
 function $h$ over a compact convex subset $K\subset \mathcal{X}
$ attains its maximum at an extreme point of $K$, i.e.,
\begin{equation*}
\sup \,h\left( K\right) =\max \,h\left( \mathcal{E}(K)\right) .
\end{equation*}%
Here, $\mathcal{E}(K)$ is the (non--empty) set of extreme points of $K$.
\end{lemma}

\noindent Indeed, by combining Proposition \ref{lemma gamma regularisation}
with Lemma \ref{Bauer maximum principle} it is straightforward to check the
following statement which does not seem to have been observed before:

\begin{lemma}[Extension of the Bauer maximum principle]
\label{Bauer maximum principle bis}\mbox{ }\newline
Let $h_{\pm }$ be two convex real--valued
 functions from a locally convex real
space $\mathcal{X}$ to $\left( -\infty ,\infty \right] $ such that $h_{-}$
and $h_{+}$ are respectively lower and upper semi--continuous. Then the
supremum of the sum $h:=h_{-}+h_{+}$ over a compact convex subset $K\subset
\mathcal{X}$ can be reduced to the (non--empty) set $\mathcal{E}(K)$ of
extreme points of $K$, i.e.,
\begin{equation*}
\sup \,h\left( K\right) =\sup \,h\left( \mathcal{E}(K)\right) .
\end{equation*}
\end{lemma}

\textit{Proof. }We first use Proposition \ref{lemma gamma regularisation} in
order to write $h_{-}=\Gamma \left( h_{-}\right) $ as a supremum over affine
and continuous functions. Then we commute this supremum with the one over $%
K$ and apply the Bauer maximum principle to obtain that%
\begin{equation*}
\sup \,h\left( K\right) =\sup \left\{ \sup \,\left[ m+h_{+}\right] (\mathcal{%
E}(K)):m\in \mathrm{A}\left( \mathcal{X}\right) \;\text{and }m|_{K}\leq
h_{-}|_{K}\right\} .
\end{equation*}%
The lemma follows by commuting again both suprema and by using $h_{-}=\Gamma
\left( h_{-}\right) $.
\hfill\qed\medbreak%

\noindent Observe, however, that under the conditions of the lemma above,
the supremum of $h=h_{-}+h_{+}$ is generally not attained on $\mathcal{E}(K)$%
.

\section{Appendix\label{Section appendix}}

\setcounter{equation}{0}%
For the reader's convenience we give here a short review on the following
subjects:

\begin{itemize}
\item Dual pairs of locally convex real spaces, see, e.g., \cite{Rudin};

\item Barycenters and $\Gamma $--regularization of real--valued functions, see,
e.g., \cite{Alfsen};

\item The Mazur and Lanford III--Robinson theorems, see \cite{LanRob,Mazur}.
\end{itemize}

\noindent These subjects are rather standard. Therefore, we keep the
exposition as short as possible and only concentrate on results used in this
paper.

\subsection{Dual Pairs of Locally Convex Real Spaces}




The notion of \emph{dual pairs} is defined as follow:

\begin{definition}[Dual pairs]
\label{dual pairs}\mbox{ }\newline
For any locally convex space $(\mathcal{X},\tau )$, let $\mathcal{X}^{\ast }$
be its dual space, i.e., the set of all continuous linear functionals on $%
\mathcal{X}$. Let $\tau ^{\ast }$ be any locally convex topology on $%
\mathcal{X}^{\ast }$. $(\mathcal{X},\mathcal{X}^{\ast })$ is called a dual
pair iff, for all $x\in \mathcal{X}$, the functional $x^{\ast }\mapsto
x^{\ast }(x)$ on $\mathcal{X}^{\ast }$ is continuous w.r.t. $\tau ^{\ast }$,
and all linear functionals which are continuous w.r.t. $\tau ^{\ast }$ have
this form.
\end{definition}

\noindent By \cite[Theorems 3.4 (b) and 3.10]{Rudin}, a typical example of a
dual pair $(\mathcal{X},\mathcal{X}^{\ast })$ is given by any locally convex
real space $\mathcal{X}$ equipped with a topology $\tau $ and $\mathcal{X}%
^{\ast }$ equipped with the $\sigma (X^{\ast },X)$--topology $\tau ^{\ast }$%
, i.e., the weak$^{\ast }$--topology. We also observe that if $(\mathcal{X},%
\mathcal{X}^{\ast })$ is a dual pair w.r.t. $\tau $ and $\tau ^{\ast }$ then
$(\mathcal{X}^{\ast },\mathcal{X})$ is a dual pair w.r.t. $\tau ^{\ast }$
and $\tau $.

\subsection{Barycenters and $\Gamma $--regularization}

The theory of compact convex subsets of a locally convex real (topological
vector) space $\mathcal{X}$ is standard. For more details, see, e.g., \cite%
{Alfsen}. An important observation is the Krein--Milman theorem (see, e.g.,
\cite[Theorems 3.4 (b) and 3.23]{Rudin}) which states that any compact
convex subset $K\subset \mathcal{X}$ is the closure of the convex hull of
the (non--empty) set $\mathcal{E}(K)$ of its extreme points. Restricted to
finite dimensions this theorem corresponds to a classical result of
Minkowski which, for any $x\in K$ in a (non--empty) compact convex subset $%
K\subset \mathcal{X}$, states the existence of a finite number of extreme
points $\hat{x}_{1},\ldots ,\hat{x}_{k}\in \mathcal{E}(K)$ and positive
numbers $\mu _{1},\ldots ,\mu _{k}\geq 0$ with $\Sigma _{j=1}^{k}\mu _{j}=1$
such that
\begin{equation}
x=\overset{k}{\sum\limits_{j=1}}\mu _{j}\hat{x}_{j}.  \label{barycenter1}
\end{equation}%
To this\ simple decomposition we can associate a probability measure, i.e.,
a \emph{normalized positive Borel regular measure}, $\mu $\ on $K$.

Borel sets of any set $K$ are elements of the $\sigma $--algebra $\mathfrak{B%
}$ generated by closed -- or open -- subsets of $K$. Positive Borel regular
measures are the positive countably additive set functions $\mu $ over $%
\mathfrak{B}$ satisfying
\begin{equation*}
\mu \left( B\right) =\sup \left\{ \mu \left( C\right) :C\subset B,\text{ }C%
\text{ closed}\right\} =\inf \left\{ \mu \left( O\right) :B\subset O,\text{ }%
O\text{ open}\right\}
\end{equation*}%
for any Borel subset $B\in \mathfrak{B}$ of $K$. If $K$ is compact then any
positive Borel regular measure $\mu $ (one--to--one) corresponds to an
element of the set $M^{+}(K)$ of Radon measures with $\mu \left( K\right)
=\left\Vert \mu \right\Vert $, and we write
\begin{equation}
\mu \left( h\right) =\int_{K}\mathrm{d}\mu (\hat{x})\;h\left( \hat{x}\right)
\label{barycenter1bis}
\end{equation}%
for any continuous function $h$ on $K$. A probability measure $\mu \in
M_{1}^{+}(K)$ is per definition a positive Borel regular measure $\mu \in
M^{+}(K)$ which is \emph{normalized}: $\left\Vert \mu \right\Vert =1$.

Therefore, using the probability measure $\mu _{x}\in M_{1}^{+}(K)$\ on $K$
defined by
\begin{equation*}
\mu _{x}=\overset{k}{\sum\limits_{j=1}}\mu _{j}\delta _{\hat{x}_{j}}
\end{equation*}%
with $\delta _{y}$ being the Dirac -- or point -- mass\footnote{$\delta _{y}$
is the Borel measure such that, for any Borel subset $B\in \mathfrak{B}$ of $%
K$, $\delta _{y}(B)=1$ if $y\in B$ and $\delta _{y}(B)=0$ if $y\notin B$.}
at $y$, Equation (\ref{barycenter1}) can be seen as an integral defined by (%
\ref{barycenter1bis}) for the probability measure $\mu _{x}\in M_{1}^{+}(K)$%
:
\begin{equation}
x=\int_{K}\mathrm{d}\mu _{x}(\hat{x})\;\hat{x}\ .  \label{barycenter2}
\end{equation}%
The point $x$ is in fact the \emph{barycenter} of the probability measure $%
\mu _{x}$. This notion is defined in the general case as follows (cf. \cite[%
Eq. (2.7) in Chapter I]{Alfsen}):

\begin{definition}[Barycenters of probability measures in convex sets]
\label{def barycenter}Let $K\subset \mathcal{X}$ be any (non--empty) compact
convex subset of a locally convex real space $\mathcal{X}$ and let $\mu \in
M_{1}^{+}(K)$ be a probability measure on $K$. We say that $x\in K$ is the
barycenter\footnote{%
Other terminologies existing in the literature: \textquotedblleft $x$ is
represented by $\mu $\textquotedblright , \textquotedblleft $x$ is the
resultant of $\mu $\textquotedblright .} of $\mu $ if, for all $z^{\ast }\in
\mathcal{X}^{\ast }$,
\begin{equation*}
z^{\ast }\left( x\right) =\int_{K}\mathrm{d}\mu (\hat{x})\;z^{\ast }\left(
\hat{x}\right) .
\end{equation*}
\end{definition}

\noindent Barycenters are well--defined for \emph{all} probability measures
in convex compact subsets of locally convex real spaces (cf. \cite[Theorems
3.4 (b) and 3.28]{Rudin}):

\begin{theorem}[Well-definiteness and uniqueness of barycenters]
\label{thm barycenter}\mbox{ }\newline
Let $K\subset \mathcal{X}$ be any (non--empty) compact subset of a locally
convex real space $\mathcal{X}$ such that $\overline{\mathrm{co}}\left(
K\right) $ is also compact. Then, for any probability measure $\mu \in
M_{1}^{+}(K)$ on $K$,$\ $there is a unique barycenter $x_{\mu }\in \overline{%
\mathrm{co}}\left( K\right) $.
\end{theorem}

\noindent Note that Barycenters can also be defined in the same way via
affine continuous functions instead of continuous linear functionals, see,
e.g., \cite[Proposition I.2.2.]{Alfsen} together with \cite[Theorem 1.12]%
{Rudin}.

It is natural to ask whether, for any $x\in K$ in the compact convex set $K$%
, there is a (possibly not unique) probability measure $\mu _{x}$ on $K$
(pseudo--) supported on $\mathcal{E}(K)$ with barycenter $x$. Equation (\ref%
{barycenter2}) already gives a first positive answer to that problem in the
finite dimensional case. The general case, which is a remarkable refinement
of the Krein--Milman theorem, has been proven by Choquet--Bishop--de Leeuw
(see, e.g., \cite[Theorem I.4.8.]{Alfsen}).

We conclude now by a crucial property concerning the $\Gamma $%
--regularization of extended real--valued
 functions in relation with the concept of
barycenters (cf. \cite[Corollary I.3.6.]{Alfsen}):

\begin{theorem}[Barycenters and $\Gamma $--regularization]
\label{Thm - Corollary I.3.6}\mbox{ }\newline
Let $K\subset \mathcal{X}$ be any (non--empty) compact convex subset of a
locally convex real space $\mathcal{X}$ and $h:K \to \mathbb{R}$
 be a continuous real--valued function. Then, for any $x\in K$,
there is a probability measure $\mu _{x}\in M_{1}^{+}(K)$\ on $K$ with
barycenter $x$ such that
\begin{equation*}
\Gamma \left( h\right) \left( x\right) =\int_{K}\mathrm{d}\mu _{x}(\hat{x}%
)\;h\left( \hat{x}\right) .
\end{equation*}
\end{theorem}

\noindent This theorem is a very important statement used to prove Theorem %
\ref{theorem trivial sympa 1}.

\subsection{The Mazur and Lanford III--Robinson Theorems}

If $\mathcal{X}$ is a separable real Banach space and $h$ is a continuous
convex real--valued
 function on $\mathcal{X}$ then it is well--known that $h$ has,
on each point $x\in \mathcal{X}$, at least one Fenchel
subgradient \textrm{d}$h\in \mathcal{X}^{\ast }$.
The Mazur theorem describes the set $\mathcal{Y}$ on
which a continuous convex function $h$ is G\^{a}teaux differentiable, more
precisely, the set $\mathcal{Y}$ for which $h$ has exactly one Fenchel
subgradient  \textrm{d}$h_{x}\in \mathcal{X}^{\ast }$ at any $x\in \mathcal{Y}$:

\begin{theorem}[Mazur]
\label{Mazur}\mbox{ }\newline
Let $\mathcal{X}$ be a separable real Banach space and let $h:\mathcal{X}%
\rightarrow \mathbb{R}$ be a continuous convex function. The set $\mathcal{%
Y}\subset \mathcal{X}$ of elements where $h$ has exactly one Fenchel
subgradient  \textrm{d}$h_{x}\in \mathcal{X}^{\ast }$
at any $x\in \mathcal{Y}$ is
residual, i.e., a countable intersection of dense open sets.
\end{theorem}

\begin{remark}
\label{Mazur remark}By Baire category theorem, the set $\mathcal{Y}$ is
dense in $\mathcal{X}$.
\end{remark}

\noindent The Lanford III--Robinson theorem \cite[Theorem 1]{LanRob}
completes the Mazur theorem by characterizing the Fenchel
subdifferential $\partial
h(x)\subset \mathcal{X}^{\ast }$ at any $x\in \mathcal{X}$:

\begin{theorem}[Lanford III -- Robinson]
\label{Land.Rob}\mbox{ }\newline
Let $\mathcal{X}$ be a separable real Banach space and let $h:\mathcal{X}%
\rightarrow \mathbb{R}$ be a continuous convex function. Then the
Fenchel
subdifferential $\partial h(x)\subset \mathcal{X}^{\ast }$ of $h$, at any $%
x\in \mathcal{X}$, is the weak$^{\ast }$--closed convex hull of the set $%
\mathcal{Z}_{x}$. Here, at fixed $x\in \mathcal{X}$, $\mathcal{Z}_{x}$ is
the set of functionals $x^{\ast }\in \mathcal{X}^{\ast }$ such that there is
a net $\{x_{i}\}_{i\in I}$ in $\mathcal{Y}$ converging to $x$ with the
property that the unique Fenchel
subgradient $\mathrm{d}h_{x_{i}}\in \mathcal{X}^{\ast }$
of $h$ at $x_{i}$ converges towards $x^{\ast }$ in the weak$^{\ast}$--topology.
\end{theorem}

\addcontentsline{toc}{section}{References}%


\begin{thebibliography}{99}


\bibitem{Zeidler3} \textsc{E. Zeidler}, \textit{Nonlinear Functional
Analysis and its Applications III: Variational Methods and Optimization.}
New York: Springer--Verlag, 1985


\bibitem{Mueller}   \textsc{S. M{\"u}ller}, \textit{Minimizing sequences for nonconvex functionals, phase transitions and singular perturbations.}
Lecture Notes in Physics Vol. 359, Springer, Englewood Cliffs,
 NJ, pp. 31--44 (1990)

\bibitem{Alfsen} \textsc{E. M. Alfsen}, \textit{Compact convex sets and
boundary integrals}. Ergebnisse der Mathematik und ihrer Grenzgebiete --
Band 57. Springer-Verlag, 1971

\bibitem{Benoist} \textsc{J. Benoist and J.-B. Hiriart-Urruty},
 \textit{What is the subdifferential of the closed
convex hull of a function?}
SIAM J. Math. Anal. Vol. 27, 1661--1679 (1996)

\bibitem{BruPedra2} \textsc{J.-B. Bru and W. de Siqueira Pedra},
Non--cooperative Equilibria of Fermi Systems With Long Range Interactions. Memoirs of the AMS \textbf{224} (2013), no. 1052.

\bibitem{Ginibre} \textsc{J. Ginibre}, On the Asymptotic Exactness of the
Bogoliubov Approximation for many Bosons Systems. Commun. Math. Phys.
\textbf{8}, 26--51 (1968)

\bibitem{Rudin} \textsc{W. Rudin}, \textit{Functional Analysis}. McGraw-Hill
Science, 1991

\bibitem{Phelps-conv} \textsc{R.R. Phelps}, \textit{Convex Functions, Monotone Operators and Differentiability}. Lecture Notes in Mathematics,
Springer-Verlag Berlin and Heidelberg, 1993

\bibitem{LanRob} \textsc{O.E. Lanford III and D.W. Robinson}, Statistical
mechanics of quantum spin systems. III. Commun. Math. Phys. \textbf{9},
327--338 (1968)

\bibitem{Simon} \textsc{B. Simon}, \textit{The Statistical Mechanics of
Lattice Gases.} Princeton: University Press, 1993

\bibitem{Mazur} \textsc{S. Mazur}, {\"{U}}ber konvexe Menge in linearen
normierten Raumen. Studia. Math. \textbf{4}, 70--84 (1933)
\end{thebibliography}
\end{document}